\documentclass[11pt,twoside]{jgsp}

\def\unu{\Id}
\def\k{\kappa}
\def\om{J}

\def\KK{\mathbb K}
\def\ZZ{\mathbb Z}

\def\RR{\mathbb R}

\def\FF{F}

\newcommand{\beq}{\begin{equation}}
\newcommand{\bee}{\end{equation}}
\newcommand{\beqa}{\begin{eqnarray}}
\newcommand{\eeqa}{\end{eqnarray}}

\newcommand{\e}{\varepsilon}

\newtheorem{theorem}{Theorem}

\def\Id{{\rm Id}\, }

\def\aut#1#2{
{\small
\noindent
\parbox[t]{6.5cm}{#1}
 \hfill
\parbox[t]{6.5cm}{#2}
}
}

\begin{document} 
 
\thispagestyle{plain} 
 
\title{Solutions for the constant quantum Yang-Baxter equation from Lie (super)algebras} 
\author{A. Tanas\u a, \'A. Ballesteros and  F. J. Herranz} 
 
\date{}

\maketitle

\comm{Communicated by XXX}\ 
 
\begin{abstract} 
We   present a systematic procedure  to obtain singular solutions of the constant quantum Yang-Baxter equation in arbitrary dimension. This approach, inspired in the Lie  (super)algebra structure,  is explicitly applied to the particular case of (graded) contractions of the orthogonal real algebra  ${\mathfrak{so}}(N+1)$.  In this way     we   show that   ``classical" contraction parameters which appear in the commutation relations of the contracted Lie algebras,  become quantum deformation parameters, arising as entries of the resulting quantum $R$-matrices.
\end{abstract}

\label{first} 
 
\section{Introduction}

Quantum $R$-matrices are solutions of the constant quantum Yang-Baxter equation (cQYBE)
 \beq
\label{ybe}
{R}_{12}{R}_{13}{R}_{23}
={R}_{23}{R}_{13}{R}_{12},
\bee
where ${R}=\sum_{i}
a_i\otimes b_i$ is a linear operator acting on a $D^2$-dimensional space and
\beq
\label{aa}
 {R}_{12}\equiv\sum_{i} a_i\otimes
b_i\otimes 1, \  {R}_{13}\equiv\sum_{i} a_i\otimes 1\otimes b_i, \  {R}_{23}\equiv\sum_{i} 1\otimes a_i\otimes b_i.
\bee
The cQYBE can be considered as a limiting case of the QYBE with spectral parameters, which constitutes the algebraic keystone for the integrability properties of (1+1) solvable models~\cite{Yang,Baxter}. Constant quantum $R$-matrices have been shown to be relevant in quantum group theory and non-commutative geometry~\cite{FRT}, since constant quantum $R$-matrices can be used to get the defining relations for non-commutative spaces such as the ones obtained under different generalizations/deformations of the special relativity theory (see~\cite{luki} and references therein). 
Several classifications for the solutions of the cQYBE, mainly concerning low dimensions, can be found in~\cite{2-1, hietarinta,2-3,2-4,2-5}. However, few constructive procedures for solutions in arbitrary dimensions $D$ are available.
The  aim of this contribution is to present a systematic construction of multiparametric solutions of the  cQYBE by means of the structure constants of any Lie (super)algebra. In Section 2 the generic $R$-matrix is constructed and  in Section 3 this approach is used to obtain explicitly the solutions generated by a family of contractions of the Lie algebra ${\mathfrak{so}}(N+1)$. In this way, by restoring to the quantum group interpretation of quantum $R$-matrices, we show that the contraction parameters (which in this case are endowed with a precise geometrical and physical meaning) can be interpreted as quantum deformation parameters in some non-commutative framework.

\section{Solutions for the constant  Quantum Yang-Baxter Equation } 

The main result of this contribution can be stated as follows.
\begin{theorem}
\label{main}
Let $X_1, \dots, X_d, X_{D}$ ($D=d+1$) span a vector space endowed with a bilinear law 
\beq
\label{comp}
X_i\ast  X_j = C_{ij}^k X_k, \  i,j,k=1,\dots, D\
\mbox{ such that } \ C_{ij}^k = 0\ \mbox{ if }\ i,j \mbox{ or } k=D,
\bee
while the remaining $C_{ij}^k$ are completely arbitrary coefficients.
Consider now the $D^2$-dimensional square    $R$-matrix with entries given by
\beq
\label{def}
R_{(i,j),(k,\ell)}= C_{ij}^k \delta_\ell^D + C_{ij}^\ell \delta_k^D,\quad
(i,j), (k,\ell) \in\{ (1,1), \dots , (D,D) \}.
\bee
Then $R$ provides a $D$-state    solution of the  cQYBE.
\end{theorem}

We stress that each non-zero coefficient $C_{ij}^k$ (\ref{comp}) is promoted into a quantum deformation one through the $R$-matrix (\ref{def}). This, in turn, means that our approach affords the construction of multiparametric $R$-matrices by simply considering different coefficients; obviously, one can take all the  
$C_{ij}^k$ equal to a single coefficient.  Furthermore, the $C_{ij}^k$ can be taken as (real or complex) constants as well as  functions depending on some other parameters, without any restriction. \\
The composition law \eqref{comp} is, in fact, a Lie (super)bracket inspired law, since the latter
can be recovered as  a particular case, once the $C_{ij}^k$  are identified with some structure constants. Thus the mechanism \eqref{def} provides a way of making a connection between a Lie (super)algebra of dimension $d$ and a $(d+1)$-state solution of the cQYBE. This is achieved    by adding a central  charge $X_D$ (an explicit application  is performed in   section 3).  Note however that the connection is only in one way, by starting from a Lie (super)algebra one can obtain a corresponding solution for the cQYBE, but the reciprocal assertion is not compulsory true.\\
By considering definition \eqref{def}, one has several lines 
(the lines $(D,j)$ and $(i,D)$, $i,j=1,\dots, D$, that is $2D-1$ lines) 
and columns (the columns $(k,\ell)$ for which both $k$ and $\ell$ are different from $D$, that is $(D-1)^2+1$ columns) which are identically zero.  Hence $\det R=0$,  whatever the  rest of the entries of the $R$-matrix with free  coefficients  $C_{ij}^k$ are, so that
 we are always dealing   with singular (non-invertible) solutions of the cQYBE

\subsection{Proof of Theorem \ref{main}}

  Let ${\cal E}^{(i,j), (k, \ell)}$ be the $D^2$-dimensional square matrix with only zero entries except for the $((i,j),  (k, \ell))$ entry, which is equal to $1$. The set $({\cal E}^{(i,j), (k, \ell)})$, $(i,j),  (k, \ell)\in\{ (1,1) \dots , (D,D) \}$)   forms a basis of the square matrices ${\cal M}_{D^2} (\KK)$ over the field $\KK$.
Note that
\beq
{\cal E}^{(i,j), (k, \ell)}= \e^{i,k}\otimes \e^{j,\ell} ,
\bee
where $\e^{i,k}$ is the $D$-dimensional square matrix with only zero entries except for the $(i,k)$ entry, which is equal to $1$. Thus $\e^{i,k}\otimes \e^{j,\ell}$, (with $i,j,  k, \ell=1,\dots, D$) is also a basis of  ${\cal M}_{D^2} (\KK)$.
Hence one can write the $R$-matrix with entries  (\ref{def}) as
\beq
R = R_{(i,j),(k,\ell)}{\cal E}^{(i,j), (k, \ell)}= R_{(i,j),(k,\ell)} \e^{i,k}\otimes \e^{j,\ell},
\bee
where hereafter we assume sum over repeated indices. Then   the 3-sites tensor product $R$-matrices (\ref{aa}), which belong to ${\cal M}_{D^3} (\KK)$,  read
\begin{align}
R_{12}&= R_{(i_1,j_1), (k_1, \ell_1)} \e^{(i_1,k_1)}\otimes \e^{(j_1,\ell_1)} \otimes \unu\notag\\
R_{13}&= R_{(i_2,j_2), (k_2, \ell_2)} \e^{(i_2,k_2)} \otimes \unu\otimes \e^{(j_2,\ell_2)}\label{ajutor}\\
R_{23}&= R_{(i_3,j_3), (k_3, \ell_3)}  \unu\otimes\e^{(i_3,k_3)} \otimes \e^{(j_3,\ell_3)},\notag
\end{align}
 where $\unu$ is the $D$-dimensional unit matrix.\\
The strategy we adopt is to   explicitly calculate the tensorial products in the
LHS and RHS of the cQYBE (\ref{ybe})  showing that both of them are identically equal to $0$. For this, one needs the following formulas (which
can be directly checked):
\begin{align}
 \e^{(i,k)}\otimes \e^{(j,\ell)}&= E^{(i-1)D+j, (k-1)D+\ell}\nonumber\\
E^{I,J}\otimes \unu &= \sum_{m=1}^D { F}^{(I-1)D+m, (J-1)D+m} \nonumber\\
\e^{(i,k)}\otimes \unu &= \sum_{m=1}^D E^{(i-1)D+m, (k-1)D+m} \label{tensor}\\
E^{I,J}\otimes \e^{j,\ell}&= F^{(I-1) D + j, (J-1) D +\ell} \nonumber\\
\unu \otimes \e^{i,k}&= \sum_{m=1}^D E^{(m-1)D+i,\, (m-1)D+k},\nonumber
\end{align}
for any $i,k, j,\ell=1,\dots, D$; $I,J=1, \dots, D^2$; and where  $E^{I,J}$ is the 
  $D^2$-dimensional square matrix with only zero entries except for the $(I,J)$ entry, which is equal to $1$, while $F^{A,B}$ ($A,B=1,\dots, D^3$) is the 
  $D^3$-dimensional square matrix with only zero entries except for the $(A,B)$ entry, which is equal to $1$. From these expressions, the matrices \eqref{ajutor} can be rewritten as 
\begin{align}
R_{12}&= R_{(i_1,j_1), (k_1, \ell_1)} \FF^{A_1, B_1}\nonumber\\
R_{13}&= R_{(i_2,j_2), (k_2, \ell_2)}\FF^{A_2, B_2} \label{123}\\
R_{23}&= R_{(i_3,j_3), (k_3, \ell_3)}\FF^{A_3, B_3} ,\nonumber\
\end{align}
where
\begin{align}
A_1&= (i_1 -1)D^2 + (j_1 -1) D + m_1 \nonumber\\
B_1&=(k_1 -1)D^2 + (\ell_1 -1) D + m_1 \nonumber\\
A_2&= (i_2 -1)D^2 + (m_2-1) D + j_2 \label{ABC} \\
B_2&= (k_2 -1)D^2 + (m_2-1) D + \ell_2 \nonumber\\
A_3&= (m_3 -1)D^2 + (i_3-1) D + j_3 \nonumber\\
B_3&= (m_3 -1)D^2 + (k_3-1) D + \ell_3.\nonumber
\end{align}
Notice that in each matrix     \eqref{123} one has five summations from    1 to $D$; four of them correspond to  the repeated indices $i_a,j_a, k_a, \ell_a$ and the other to   $m_a$  ($a=1,2,3$).\\
By taking into account the   property
\beqa
\label{ultima}
\FF^{A_1,B_1}\FF^{A_2, B_2}=\delta^{A_2, B_1} \FF^{A_1, B_2},
\eeqa
we obtain that   the LHS of the cQYBE \eqref{ybe} contains the matrix product 
\beqa
\FF^{A_1, B_1}\FF^{A_2, B_2}\FF^{A_3, B_3}=\delta^{A_2, B_1} \delta^{A_3, B_2} \FF^{A_1, B_3}.
\label{parte}
\eeqa
Since   the dimension $D$ is arbitrary, equations \eqref{ABC} imply that
\beq
\begin{aligned}
\delta^{A_2,B_1}&=\delta^{k_1, i_2}\delta^{\ell_1, m_2} \delta^{j_2, m_1}\\
\delta^{A_3,B_2}&=\delta^{k_2, m_3}\delta^{i_3, m_2} \delta^{\ell_2, j_3}.
\end{aligned}
\label{delte}
\bee
By firstly inserting  (\ref{parte}),  (\ref{delte})  and  the matrix elements $R_{(i,j),(k,\ell)}$ \eqref{def} in the LHS of the cQYBE \eqref{ybe}, and secondly expanding the terms one finally obtains that
\beqa
\label{LHS}
R_{12}R_{13}R_{23}=0,
\eeqa
since $C_{D j_a}^{k_a}=C_{i_a D}^{k_a}=C_{i_a j_a}^{D}=0$.
Similarly, one also find that
$R_{23}R_{13}R_{12}=0$, so that  the $R$-matrix defined by theorem \ref{main} is a $D$-state solution of the cQYBE. QED

\section{$R$-matrices from Contractions of  Orthogonal Lie Algebras} 

Theorem \ref{main} provides some general results   to construct singular solutions of the cQYBE. Nevertheless, as commented above, this can be specially applied to Lie (super)algebras by simply introducing their structure constants in the definition of the $R$-matrix entries (\ref{def}). In this Section we construct explicitly the $R$-matrices corresponding to a  particular family of   contracted algebras which are obtained from  ${\mathfrak{so}}(N+1)$. As a byproduct, we find that the graded contraction parameters are promoted into quantum deformation ones within such $R$-matrices.\\
Let us consider the real Lie algebra  ${\mathfrak{so}}(N+1)$ whose $\frac 12 N(N+1)$ generators
$\om_{ab}$ $(a,b=0,1,\dots, N$,
$a<b)$ satisfy the   non-vanishing Lie brackets given by  
\beq
 [\om_{ab}, \om_{ac}] =  \om_{bc}  ,\quad [\om_{ab}, \om_{bc}] = -\om_{ac}  ,\quad [\om_{ac}, \om_{bc}]
=  \om_{ab}  ,
\label{aaa}  
\bee 
where $a<b<c$.
The   $\ZZ_2^{\otimes N}$-graded contractions of ${\mathfrak{so}}(N+1)$ contain  the so called Cayley-Klein (CK) orthogonal  Lie algebras~\cite{marc}. This family, denoted collectively      ${\mathfrak{so}}_{\k}(N+1)$,        depends   on
$N$ real contraction parameters $\k=(\k_1,\dots,\k_N)$. The non-zero commutators   turn out to be~\cite{marc}: 
\beq 
[\om_{ab}, \om_{ac}] =  \k_{ab}\om_{bc} , \quad [\om_{ab}, \om_{bc}] = -\om_{ac}  ,\quad [\om_{ac},
\om_{bc}] =  \k_{bc}\om_{ab},  \label{aab}  
\bee 
 without sum over repeated indices, and where  the two-index parameters $\k_{ab}$ are expressed in
terms of the
$N$ basic  ones through
\beq
\k_{ab}=\k_{a+1}\k_{a+2}\cdots\k_b ,\quad a,b=0,1,\dots,N, \quad a<b.
\label{aac}
\bee
Each   contraction parameter $\k_\mu$ can take a positive, negative or zero value, so that ${\mathfrak{so}}_{\k}(N+1)$  comprises $3^N$ Lie algebras (some of them are isomorphic).
For instance~\cite{casimir},   when
$\k_\mu\ne 0$ $\forall \mu$,  ${\mathfrak{so}}_{\k}(N+1)$   is a
 simple pseudo-orthogonal  algebra ${\mathfrak{so}}(p,q)$ ($p+q=N+1$)     (the $B_l$ and $D_l$ Cartan series);  when    $\k_1=0$ we recover the    inhomogeneous algebras 
 ${\mathfrak{iso}}(p',q')$ ($p'+q'=N$); and when all   $\k_\mu=0$, we find the flag algebra $ {\mathfrak{ i\dots iso}}(1)$.
We recall that kinematical algebras, such us Poincar\'e, Galilei, (anti-)de Sitter, etc., associated to different models of spacetimes of constant
curvature  also belong to this CK family of algebras~\cite{casimir}. \\
Now, in order to apply Theorem \ref{main}  we enlarge the CK algebra ${\mathfrak{so}}_{\k}(N+1)$ with an additional central generator $\Xi$, that is  $[\Xi,J_{ab}]=0$, $\forall ab$. The vector space 
corresponding  to ${\mathfrak{so}}_{\k}(N+1)\oplus \RR$ is spanned by $D=\frac 12 N(N+1)+1$ elements.  We label the $D$ generic generators  $\{X_1,,\dots, X_{\frac 12 N(N+1)}, X_{D}\}$ 
as $\{J_{01},\dots,J_{N-1\,N},\Xi\equiv X_D\}$ according to the increasing    order of $ab$ with $a<b$;
the indices of the entries   (\ref{def}) run as 
 $i,j,k,\ell=\{ 01,02,\dots,0N,12,\dots,N-1\,N,D\}$.  Then by taking into account the  structure constants of (\ref{aab}) we obtain that the $D$-state solution of the cQYBE associated to ${\mathfrak{so}}_{\k}(N+1)$  is the $R$-matrix with the following non-zero entries:
\beq
\begin{array}{ll}
R_{(ab, ac),(bc,D)}= \k_{ab},&\quad R_{(ab, ac),(D ,bc)}= \k_{ab}\\
R_{(ac ,ab),(bc,D)}= -\k_{ab},&\quad R_{(ac ,ab),(D ,bc)}= -\k_{ab}\\
R_{(ab, bc),(ac,D)}= -1,&\quad R_{(ab, bc),(D, ac)}= -1 \\
R_{(bc ,ab),(ac,D)}= 1,&\quad R_{( bc, ab),(D, ac)}= 1 \\
R_{(ac, bc),(ab,D)}= \k_{bc},&\quad R_{(ac, bc),(D , ab)}= \k_{bc}\\
R_{(bc ,ac),(ab,D)}=- \k_{bc},&\quad R_{(bc ,ac),(D , ab)}= -\k_{bc},
 \end{array}
\label{aad}
\bee
where $a,b,c=0,1,\dots,N$ and $a<b<c$.
Thus we have obtained a multiparametric solution of the cQYBE, which holds simultaneously for the $3^N$ particular Lie algebras contained in the CK family. The maximum number of quantum deformation parameters is   $N$  $(\k_1,\dots,\k_N)$, that corresponds to   ${\mathfrak{so}}(p,q)$; through the  contractions $\k_\mu=0$ this number is subsequently reduced up to reach the flag algebra, for which there is no quantum parameters other than the constants $\pm 1$.\\
Let us illustrate explicitly this construction with the  $N=2$ case.

\subsection{Solutions of the cQYBE associated to ${\mathfrak{so}}_{\k_1,\k_2}(3)$ } 

The CK algebra with $N=2$ is ${\mathfrak{so}}_{\k_1,\k_2}(3)$, which depends on two real coefficients $(\k_1,\k_2)$ and is spanned   by  three generators $\{ J_{01}, J_{02},J_{12}\}$ fulfilling 
\beq
\label{so3}
\left[ J_{01}, J_{02}\right] = \k_1 J_{12}, \quad 
\left[ J_{01}, J_{12}\right] =- J_{02}, \quad
\left[ J_{02}, J_{12}\right] = \k_2 J_{01}.
\bee
According to  the pair $(\k_1,\k_2)$ we find that ${\mathfrak{so}}_{\k_1,\k_2}(3)$ covers     9   Lie algebras:
${\mathfrak{so}}(3)$ for $(+,+)$; ${\mathfrak{so}}(2,1)$ for $(+,-)$, $(-,+)$ and $(-,-)$;     ${\mathfrak{iso}}(2)$ for $(+,0)$ and $(0,+)$; ${\mathfrak{iso}}(1,1)$ for $(-,0)$ and $(0,-)$; and   
${\mathfrak{iiso}}(1)$ for $(0,0)$. By considering the Lie group ${SO}_{\k_1,\k_2}(3)$,    both contraction parameters, $\k_1,\k_2$, can be identified with the constant  curvature of the two-dimensional homogeneous space of points ${SO}_{\k_1,\k_2}(3)/\langle J_{12}\rangle$ and of lines ${SO}_{\k_1,\k_2}(3)/\langle J_{01}\rangle$,  respectively~\cite{ck2}.
 Furthermore, if  $\{ J_{01}, J_{02},J_{12}\}$ are interpreted, in this order, as time-translation, spatial-translation and boost generators, then the six CK algebras with $\k_2\le 0$ are kinematical ones~\cite{ck2}. In this case, besides the geometrical interpretation, the contraction parameters have a physical meaning as well, since they can be expressed as $\k_1=\pm 1/\tau^2$ where $\tau$ is the universe radius and $\k_2=-1/c^2$ where $c$ is the speed of light. Thus 
 ${\mathfrak{so}}_{\k_1,\k_2}(3)$ and ${SO}_{\k_1,\k_2}(3)/\langle J_{12}\rangle$  comprises the following kinematical algebras and (1+1) spacetimes~\cite{ck2}: the anti-de Sitter $(+1/\tau^2,-1/c^2)$, Minkowskian $(0,-1/c^2)$, de Sitter $(-1/\tau^2,-1/c^2)$,   oscillating Newton-Hooke $(+1/\tau^2,0)$, Galilean $(0,0)$ and 
 expanding Newton-Hooke ones $(-1/\tau^2,0)$. \\
 Next we present the $D=4$-state solution of the cQYBE coming from ${\mathfrak{so}}_{\k_1,\k_2}(3)$.
 At this dimension the indices $i,j,k,\ell=\{ 01,02,12,D=4\}$, so that the  entries (\ref{aad}) give rise 
to  the following $16\times 16$ $R$-matrix, $R=$
\[ \left( \begin{array}{cccc|cccc|cccc|cccc}
0 & 0 & 0 & 0 & 0 & 0 & 0 & 0 & 0 & 0 & 0 & 0 & 0 & 0 & 0 & 0\\
0 & 0 & 0 & 0 & 0 & 0 & 0 & 0 & 0 & 0 & 0 & \k_1 & 0 & 0 & \k_1 & 0\\
0 & 0 & 0 & 0 & 0 & 0 & 0 & -1 & 0 & 0 & 0 & 0 & 0 & -1 & 0 & 0\\
0 & 0 & 0 & 0 & 0 & 0 & 0 & 0 & 0 & 0 & 0 & 0 & 0 & 0 & 0 & 0\\
\hline
0 & 0 & 0 & 0 & 0 & 0 & 0 & 0 & 0 & 0 & 0 & -\k_1 & 0 & 0 & -\k_1 & 0\\
0 & 0 & 0 & 0 & 0 & 0 & 0 & 0 & 0 & 0 & 0 & 0 & 0 & 0 & 0 & 0\\
0 & 0 & 0 & \k_2 & 0 & 0 & 0 & 0 & 0 & 0 & 0 & 0 & \k_2 & 0 & 0 & 0\\
0 & 0 & 0 & 0 & 0 & 0 & 0 & 0 & 0 & 0 & 0 & 0 & 0 & 0 & 0 & 0\\
\hline
0 & 0 & 0 & 0 & 0 & 0 & 0 & 1 & 0 & 0 & 0 & 0 & 0 & 1 & 0 & 0\\
0 & 0 & 0 & -\k_2 & 0 & 0 & 0 & 0 & 0 & 0 & 0 & 0 & -\k_2 & 0 & 0 & 0\\
0 & 0 & 0 & 0 & 0 & 0 & 0 & 0 & 0 & 0 & 0 & 0 & 0 & 0 & 0 & 0\\
0 & 0 & 0 & 0 & 0 & 0 & 0 & 0 & 0 & 0 & 0 & 0 & 0 & 0 & 0 & 0\\
\hline
0 & 0 & 0 & 0 & 0 & 0 & 0 & 0 & 0 & 0 & 0 & 0 & 0 & 0 & 0 & 0\\
0 & 0 & 0 & 0 & 0 & 0 & 0 & 0 & 0 & 0 & 0 & 0 & 0 & 0 & 0 & 0\\
0 & 0 & 0 & 0 & 0 & 0 & 0 & 0 & 0 & 0 & 0 & 0 & 0 & 0 & 0 & 0\\
0 & 0 & 0 & 0 & 0 & 0 & 0 & 0 & 0 & 0 & 0 & 0 & 0 & 0 & 0 & 0
\end{array} \right). \]
Now, if we consider this $R$-matrix as the structure constant matrix for a non-commutative space constructed by using the standard FRT approach~\cite{FRT}, we get  a direct relationship between "classical" contraction/curvature parameters and quantum deformation ones. Moreover, we find that physical classical quantities such as $\tau$ and $c$ can also be promoted into quantum deformation parameters.
The  construction  of such quantum spaces is currently in progress.

\section*{Acknowledgments}

This work was partially supported  by the Spanish MEC   (FIS2004-07913) and JCyL   (VA013C05).
 A. Tanas\u{a} would like to thank the staff of the Physics Department of the University of Burgos for their hospitality during his stay.

\aut{ Adrian Tanas\u a\\ 
Labo. Physique Th\'eorique\\
b\^at. 210, CNRS UMR 8627\\
Universit\'e Paris XI\\
91405 Orsay Cedex, France\\
{\it E-mail address}:\\ 
{\tt adrian.tanasa@ens-lyon.org }
} 
{  \'Angel Ballesteros and Francisco J. Herranz\\ 
Depto. de Fisica\\
Universidad de Burgos\\
09001 Burgos, Spain\\
{\it E-mail address}:
 {\tt angelb@ubu.es}\\
 {\tt fjherranz@ubu.es}\\
 }

\label{last} 
\end{document}